\documentclass[a4paper]{jpconf}
\usepackage{amsfonts, amsmath, amssymb, amsthm, amsbsy, amscd, euscript}
\usepackage{float}
\usepackage{pifont} % for ding{51] and ding{55}
\usepackage{textcomp} % for \textnumero

\theoremstyle{definition}
\newtheorem{example}{Example}

\theoremstyle{remark}
\newtheorem{rem}{Remark}

\newcommand{\pinner}{\mathbin{\mathchoice
   {\hbox{\vrule width0.6em depth0pt height0.4pt
   \vrule width0.4pt depth0pt height0.8ex}}
   {\hbox{\vrule width0.6em depth0pt height0.4pt
   \vrule width0.4pt depth0pt height0.8ex}}
   {\hbox{\kern0.14em
   \vrule width0.48em depth0pt height0.4pt
   \vrule width0.4pt depth0pt height0.6ex\kern0.14em}}
   {\hbox{\kern0.1em
   \vrule width0.39em depth0pt height0.4pt
   \vrule width0.4pt depth0pt height0.5ex\kern0.1em}}}}
\newcommand{\inner}{\pinner\,}

\newcommand{\veps}{\varepsilon}
\newcommand{\BBR}{\mathbb{R}}

\newcommand{\bx}{\boldsymbol{x}}

\newcommand{\dd}{\partial}

\newcommand{\Id}{{\mathrm d}}

\newcommand{\cP}{\mathcal{P}}
\newcommand{\cQ}{\mathcal{Q}}
\newcommand{\cX}{{\EuScript X}}%{\mathcal{X}}
\newcommand{\bcP}{{\boldsymbol{\mathcal{P}}}}

\newcommand{\schouten}[1]{\lshad {#1} \rshad}
\newcommand{\rP}{{\mathrm P}}
\newcommand{\vmark}{\ding{51}}
\newcommand{\xmark}{\ding{55}}

\begin{document}

%%%%%%%%%%% REMOVE IT BEFORE SUBMISSION %%%%%%%%%%%%%%%%%%%%%%%%%%%
\pagestyle{plain}
%%%%%%%%%%%%%%%%%%%%%%%%%%%%%%%%%%%%%%%%%%%%%%%%%%%%%%%%%%%%%%%%%%%

\title%[The $1:6$ ratio tetrahedral flow preserves the space of Poisson bi%\/-\/vectors]
{Do the Kontsevich tetrahedral flows preserve %\\[2pt] 
or destroy the space of Poisson bi\/-\/vectors\,?}

\author{Anass Bouisaghouane and Arthemy V Kiselev}

\address{Johann Bernoulli Institute for Mathematics and Computer Science,
  University of Groningen,
  P.O.Box~407, 9700\,AK Groningen, The Netherlands}

\ead{A.V.Kiselev@rug.nl}

\begin{abstract}\\
%\parbox[c][1.2\height][c]{13.4cm}{% The original value was 12cm, now it fits the text of the abstract in width.
\noindent%\mbox
{\rm{From the paper ``Formality Conjecture'' (Ascona 1996)%
%$($see Ref.~\cite{Ascona96}$)$}
:}\\[1pt]
\mbox{ }\quad
\textit{I am aware of only one such a class, it corresponds to simplest good graph, the complete graph with $4$~vertices $($and $6$~edges$)$. 
This class gives a remarkable vector field on the space of bi-vector fields 
on~$\mathbb{R}^{d}$. The evolution with respect to the time~$t$ is described by the following non\/-\/linear partial differential equation}{\rm:} 
$\ldots${\rm,} %{\rm$[$see %~formula %\eqref{EqFlow1} 
%below$]$,} 
\textit{where $\alpha=\sum_{i,j}\alpha_{ij}%\tfrac
{\partial}/{\partial x_{i}}\wedge %\tfrac
{\partial}/{\partial x_{j}}$ is a bi\/-\/vector field on~$\mathbb{R}^d$.%~$\ldots$
}\\ 
\mbox{ }\quad
\textit{It follows from general properties of cohomology that $1)$ \textbf{this evolution preserves the class of $\mathbf{(}$real\/-\/analytic$\mathbf{)}$ Poisson structures},~$\ldots$}\\
\mbox{ }\quad
\textit{In fact, I cheated a little bit. In the formula for the vector field on the space of bivector fields which one get from the tetrahedron graph, an additional term is present. 
$\ldots$ %This term is equal $($up to a numerical factor$)$ to
%{\rm$[$see~\eqref{EqFlow2} below$]$}. 
It is possible to prove formally that \textbf{if $\alpha$ is a Poisson %bi-vector
bracket, i.e.\ if $[{\alpha,\alpha}]=0\in T^2(\mathbb{R}^d)$, then the additional term shown above vanishes.}}%\bold
}%\end{minipage}

\smallskip\noindent%
By using twelve Poisson structures with high\/-\/degree polynomial coefficients
as explicit coun\-ter\-ex\-am\-p\-les, we show that both the above claims are false: neither does the first flow pre\-serve the property of bi\/-\/vectors to be Poisson nor does the second flow vanish identically at %the 
Poisson bi\/-\/vectors. The counterexamples at hand %themselves 
suggest a correction to the for\-mula for the ``exotic" flow on the space of Poisson bi\/-\/vectors; in fact, this flow is encoded by the balanced sum involving both the Kontsevich tetrahedral graphs (that give rise to the flows men\-ti\-o\-ned above). We reveal %prove %argue
that it is only the balance $1:6$ for which the flow does preserve the space of Poisson bi\/-\/vectors. 
\end{abstract}

%\subsection*
\paragraph{\textbf{\textup{Introduction.}}}
The Kontsevich graph complex is the language of deformation quantisation on finite\/-\/dimensional Poisson manifolds~\cite{Ascona96,KontsevichFormality}. %Let us
We consider the class of oriented graphs on %with 
two sinks and $k\geqslant1$ internal vertices (of which, each is the tail of two edges and carries a copy of the Poisson bi\/-\/vector~$\cP$). Encoding bi\/-\/differential operators, such graphs determine the flows on the space of bi\/-\/vectors on a %the 
Poisson manifold at hand. The two flows with $k=4$ %the minimal number of 
internal vertices in the graphs are provided by the two tetrahedra~\cite{Ascona96},
see Fig.~\ref{FigTetra} on the next page. %below.
%\noindent
By producing $12$ counterexamples, %of Poisson bi-vectors with polynomial coefficients of degree $\geq3$,
we prove that the claim~\cite{Ascona96,KontsevichFormality} 
of preservation of the Poisson property is false as stated.
%so that the (variational) Poisson bi\/-\/vectors are fragile with respect to the Kontsevich tetrahedral flows.
Simultaneously, we reveal that the flow which is determined by the second graph is not always %everywhere\ignore{instead of always/necessarily} 
vanishing by virtue of the skew\/-\/symmetry and Jacobi identity for %the 
Poisson bi\/-\/vectors~$\cP$.

This paper is structured as follows. First we recall the correspondence between graphs and polydifferential operators~\cite{MKZurichICM,MKParisECM} and we indicate the mechanism for such an operator to vanish% 
%by virtue of skew-symmetry and the Jacobi identity for a given Poisson bi-vector (which is contained in every internal vertex of graphs
, cf.~\cite{sqs15,f16}. 
In section~\ref{generators} we recall three constructions of %the 
Poisson brackets with polynomial coefficients of arbitrarily high degree 
(see~\cite{Donin,%Nutku
Perelomov,Van}). %, \cite{Nambu}
%%% Next, in section~\ref{SecDeform}
%%% we recite the basics of %the 
%%% Poisson structure deformation~\cite{Gerstenhaber}. 
In Tables~\ref{DonNutTab}\/--\/%, \ref{DonTab}, \ref{VanTab} and 
\ref{HDTab} on pp.~\pageref{DonNutTab}\/--\/\pageref{HDTab}
we then summarise the properties of all structures from our 
$12$~counterexamples to the claim \cite{Ascona96} that 
\begin{itemize}%{enumerate}
\item[(\textit{i})] the flow $\dot{\cP}=\Gamma_{1}(\cP)$ which the first graph in Fig.~\ref{FigTetra} encodes on the space of bi\/-\/vec\-tors~$\cP$ would preserve their property to be Poisson (in fact, it does not), and that
\item[(\textit{ii})] the flow $\dot{\cP}=\Gamma_{2}(\cP)$ would always be trivial whenever the bi\/-\/vector~$\cP$ is Poisson (in fact, this is not true). 
\end{itemize}%{enumerate}
In particular, the twelfth counterexample pertains to the infinite\/-\/dimensional jet\/-\/space geometry of variational Poisson structures~\cite{Olver}. (Quoted from \cite{Kent}, the Hamiltonian dif\-fe\-ren\-ti\-al operator for that variational Poisson bi\/-\/vector~$\bcP$ is %then 
processed by using the techniques from~\cite{gvbv,dq15,prg15}). %, cf.~\cite{prg15}.)
%In the extensive appendix\marginpar{Edit} 
%on pp.{pApp}-{p.AppEnd} we list the formulas for the bi-vectors and their respective flows. 

Finally, we examine at which balance the linear combination of the Kontsevich tetrahedral flows %whether the respective flows 
preserves the space %property 
of Poisson %bi-vectors to be . 
structures on finite\/-\/dimensional manifolds. We argue %prove 
that the ratio $1:6$ does the~job; this claim has been proved in~\cite{f16}.

\section{The graphs and operators}\label{SecGraphOp}\noindent%
Let us formalise a way to encode polydifferential operators using oriented graphs. Con\-si\-der the space $\mathbb{R}^{n}$ with Cartesian coordinates $\bx=(x_{1}$,\ $\ldots$,\ $x_{n})$, here $2\leqslant n < \infty$; for typographical reasons %only do 
we use the lower indices to enumerate the variables, so that $x_{1}^{2}=(x_{1})^{2}$, etc. By definition, the decorated %indexed 
edge $\bullet \xrightarrow{\ i\ } \bullet$ denotes at once the derivation ${\partial}/{\partial x_{i}}\equiv \partial_{i}$ (that acts on the content of the arrowhead vertex) and the summation $\sum_{i=1}^{n}$ (over the index~$i$ in %within 
the object which is contained in %within 
the ar\-row\-tail vertex). For example, the graph $\bullet \xleftarrow{\ i\ } \cP^{ij}(\bx) \xrightarrow{\ j\ } \bullet$ encodes the bi\/-\/differential operator $\sum_{i=1}^{n}(\cdot) \overleftarrow{\partial_{i}}\, \cP^{ij}(\bx)\, \overrightarrow{\partial_{j}}(\cdot)$. If its coefficients~$\cP^{ij}$ are antisymmetric, then the graph $\bullet \xleftarrow{\ i\ } \bullet \xrightarrow{\ j\ } \bullet$ encodes the 
bi\/-\/vector $\cP=\cP^{ij}\,\partial_{i}\wedge \partial_{j}$, where $\partial_{i}\wedge \partial_{j}=\tfrac{1}{2}(\partial_{i} \otimes \partial_{j}-\partial_{j}\otimes \partial_{i})$. It then specifies the Poisson bracket $\{{\cdot},{\cdot}\}_{\cP}$ if the $\tfrac{n(n-1)}{2}$-\/tuple of coefficients solves the system of equations
\begin{equation}\label{EqJac}
(\cP^{ij})\overleftarrow{\partial_{\ell}}\cdot \cP^{\ell k}
+(\cP^{jk})\overleftarrow{\partial_{\ell}}\cdot \cP^{\ell i}
+(\cP^{ki})\overleftarrow{\partial_{\ell}}\cdot \cP^{\ell j}=0, 
\end{equation}
%that is, 
hence the bracket $%\smash
{\bullet \xleftarrow[L]{\ i\ } \cP^{ij} \xrightarrow[R]{\ j\ } \bullet}$ satisfies the Jacobi identity.  
%Or equivalently, the bi-vector $\cP=\frac{1}{2} \langle \xi_{i},\cP^{ij}\xi_{j} \rangle$ satisfies the classical master-equation $\schouten{\cP,\cP}=0$ (here $\xi_{i}$ is the parity-odd canonical conjugate of $x_{i}$ for every $i$, and $\schouten{\cdot,\cdot}$ is the parity odd Schouten bracket, see \cite{gvbv\ignore{cycle14}} for discussion). 
Clearly, %we then have 
$\cP^{ij}(\bx)=\{x_{i},x_{j}\}_{\cP}$.

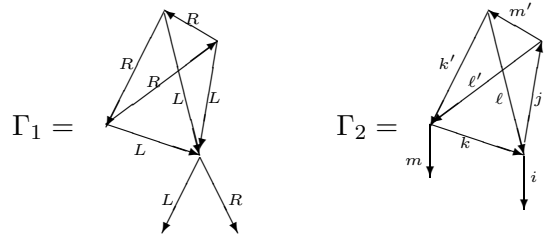
\begin{figure}[b]
\centering
\unitlength=1mm
\special{em:linewidth 0.4pt}
\linethickness{0.4pt}
\begin{picture}(30.00,27.00)(0,3)
\put(20.00,30.00){\vector(-1,-2){7.67}}
\put(12.33,15.00){\vector(4,3){14.67}}
\put(20.00,30.1){\vector(-4,3){0}}
\put(27.00,26){\line(-5,3){7.00}}
\put(20.00,29.67){\vector(1,-4){4.67}}
\put(24.7,10.65){\vector(1,-2){5.00}}
\put(24.7,10.65){\vector(-1,-2){5.00}}
\put(12.33,15.00){\vector(3,-1){12.33}}
\put(27.00,26.00){\line(-1,-6){2.2}}
\put(24.85,12.5){\vector(-1,-4){0.2}}
\put(0.00,14.67){\makebox(0,0)[lc]{$\Gamma_1={}$}}
\put(16,22.33){\makebox(0,0)[rb]{\tiny$R$}}
\put(22.7,17.67){\makebox(0,0)[rb]{\tiny$L$}}
\put(22.67,28.33){\makebox(0,0)[lb]{\tiny$R$}}
\put(25.67,18.5){\makebox(0,0)[lc]{\tiny$L$}}
\put(19.5,20.00){\makebox(0,0)[rb]{\tiny$R$}}
\put(17.67,12.33){\makebox(0,0)[rt]{\tiny$L$}}
\put(21.33,5.00){\makebox(0,0)[rc]{\tiny$L$}}
\put(28.33,5.00){\makebox(0,0)[lc]{\tiny$R$}}
\end{picture}
\qquad
\unitlength=1mm
\special{em:linewidth 0.4pt}
\linethickness{0.4pt}
\begin{picture}(30.00,27.00)(0,3)
\put(20.00,30.00){\vector(-1,-2){7.67}}
\put(27,26){\vector(-4,-3){14.5}}
\put(20.00,30.1){\vector(-4,3){0}}
\put(27.00,26){\line(-5,3){7.00}}
%\put(25.00,24.67){\vector(-1,1){5.00}}
\put(20.00,29.67){\vector(1,-4){4.67}}
\put(24.7,10.65){\vector(0,-1){7.00}}
\put(12.33,15){\vector(0,-1){7.00}}
\put(12.33,15.00){\vector(3,-1){12.33}}
\put(27.00,26.00){\line(-1,-6){2.22}}
%\put(25.00,11){\line(1,6){2}}
\put(27.00,26){\vector(1,4){0}}
\put(0.00,14.67){\makebox(0,0)[lc]{$\Gamma_2={}$}}
\put(15.83,22.33){\makebox(0,0)[rb]{\tiny$k'$}}
\put(22.00,17.67){\makebox(0,0)[rb]{\tiny$\ell$}}
\put(22.67,29.33){\makebox(0,0)[lb]{\tiny$m'$}}
\put(26,18.33){\makebox(0,0)[lc]{\tiny$j$}}
\put(19.33,20.00){\makebox(0,0)[rb]{\tiny$\ell'$}}
\put(17.67,13.33){\makebox(0,0)[rt]{\tiny$k$}}
\put(11.33,10.00){\makebox(0,0)[rc]{\tiny$m$}}
\put(25.5,8.00){\makebox(0,0)[lc]{\tiny$i$}}
\end{picture}
\caption{These tetraheral graphs encode flows~\eqref{EqFlow1} and~\eqref{EqFlow2}, respectively. Each oriented edge carries a summation index that runs from~$1$ to the dimension of the Poisson manifold at hand. For each internal vertex (where a copy of the Poisson bi\/-\/vector~$\cP$ is stored), the pair of out\/-\/going edges is ordered, L${}\prec{}$R: the left edge~(L) carries the first index and the other edge~(R) carries the second index in the bi\/-\/vector co\-ef\-fi\-ci\-ents. %, see section~\ref{SecGraphOp}. 
(In retrospect, the ordering and labelling of the in\-de\-x\-ed oriented edges can be guessed from formulas~\eqref{EqFlows} on p.~\pageref{EqFlows}.)}\label{FigTetra}
\end{figure}

From now on, let us consider only the oriented graphs whose vertices are either sinks, with no issued edges, or tails for an ordered pair of arrows, 
each decorated with %carrying 
its own index (see Fig.~\ref{FigTetra}). %on p.~\pageref{FigTetra}). 
% Figure of tetra's
%%%
Allowing the only exception in footnote~\ref{FootSkew}, 
%on p.~\pageref{FootSkew} below, 
we shall always assume that there are neither tadpoles, nor double oriented edges, nor two\/-\/edge loops.
%so that none of the three graphs which are shown here (or similar graphs) will be considered in what follows.
% Missing figure of 3 bad graphs

We also postulate that every vertex which is not a sink carries a copy of a given Poisson bi\/-\/vector $\cP=\cP^{ij}(\bx)\,\partial_{i} \wedge \partial_{j}$; the ordering of decorated %indexed 
out\/-\/going edges coincides with the ordering
``first${}\prec{}$se\-c\-ond'' of the indexes in the coefficients of~$\cP$.

\begin{example}\label{ExFlows}
Under all these assumptions, the two tetrahedra which are portrayed 
in Fig.~\ref{FigTetra} are, up to a symmetry, the only admissible graphs 
with $k=4$~internal vertices, $2k=6+2$~edges, and two sinks.
The first graph in Fig.~\ref{FigTetra} encodes the bi\/-\/vector
\begin{subequations}\label{EqFlows}
\begin{align}\label{EqFlow1}
\Gamma_{1}(\cP)&=\sum_{i,j=1}^{n} %\left
\biggl(\sum\limits_{k,\ell,m,k',\ell',m'=1}^{n} \frac{\partial^{3}\cP^{ij}}{\partial x_{k} \partial x_{\ell}  \partial x_{m}} \frac{\partial \cP^{kk'} }{\partial x_{\ell'}} \frac{\partial \cP^{\ell\ell'}}{\partial x_{m'}} \frac{\partial \cP^{mm'}}{\partial x_{k'}}%\right
\biggr)\frac{\partial}{\partial x_{i}} \wedge \frac{\partial}{\partial x_{j}}.\\
\intertext{Likewise, the second graph in Fig.~\ref{FigTetra} yields the 
bi\/-\/vector}
\label{EqFlow2}
\Gamma_{2}(\cP)&=\sum_{i,m=1}^{n} %\left
\biggl( \sum\limits_{j,k,\ell,k',\ell',m'=1}^{n} \frac{\partial^{2}\cP^{ij}}{\partial x_{k} \partial x_{\ell}} \frac{\partial^{2}\cP^{km}}{\partial x_{k'} \partial x_{\ell'} }  \frac{\partial\cP^{k'\ell}}{\partial x_{m'}}   \frac{\partial\cP^{m'\ell'}}{\partial x_{j}}   %\right
\biggr)\frac{\partial}{\partial x_{i}} \wedge \frac{\partial}{\partial x_{m}}.
\end{align}
\end{subequations}
In this paper we examine
\begin{itemize}%{enumerate}
\item[(\textit{i})] whether the respective flows $\frac{\Id}{\Id\varepsilon}(\cP)=\Gamma_{\alpha}(\cP)$ at $\alpha=1,2$ preserve or, in fact, destroy the property of bi\/-\/vectors $\cP(\varepsilon)$ to be Poisson, provided that the Cauchy datum $\cP{\bigl|}_{\varepsilon=0}$ is such;
\item[(\textit{ii})] we also inspect
whether the second flow is (actually, it is not) vanishing identically at all $\varepsilon$, provided that the Cauchy datum is a Poisson bi\/-\/vector.
\end{itemize}%{enumerate}
\end{example}

% Incomplete remark, ice cream cone graph is missing

\begin{rem}
Whenever the bi\/-\/vector~$\cP$ in every internal vertex of a non\/-\/empty graph~$\Gamma$ is Poisson, the bi\/-\/differential operator which is encoded by $\Gamma$ can vanish identically. First, this occurs due to the skew\/-\/symmetry of coefficients of the bi\/-\/vector.%\footnote{For example, consider the oriented graph with ordered pairs of indexed edges $(i<j, k<l, m<n, p<q)$: %This is hell to work out \label{FootSkew}}
\footnote{\label{FootSkew}%
For example, consider this oriented graph with ordered pairs of indexed edges\\[1pt]
%\noindent
\begin{minipage}{0.712\textwidth}
($i\prec j$,\ $k\prec\ell$,\ $m\prec n$,\ $p\prec q$). 
We claim that due to the antisymmetry of~$\cP$ which is contained in each of the four internal vertices, the operator (which this graph encodes) vanishes identically. Indeed, it equals minus itself:%\begin{subequations}
\begin{align*}
\partial_{m}&\partial_{n}(\cP^{pq})\partial_{p}(\cP^{km})\partial_{q}(\cP^{\ell n})\partial_{k}\partial_{\ell}(\cP^{ij})\,\partial_{i}\wedge \partial_{j}
\\
{}&=-\partial_{m}\partial_{n}(\cP^{qp})\partial_{p}(\cP^{km})\partial_{q}(\cP^{\ell n})\partial_{k}\partial_{\ell}(\cP^{ij})\,\partial_{i}\wedge \partial_{j}%={}
\\
{}&=-\partial_{n}\partial_{m}(\cP^{pq})\partial_{q}(\cP^{\ell n})\partial_{p}(\cP^{km})\partial_{\ell}\partial_{k}(\cP^{ij})\,\partial_{i}\wedge \partial_{j}=0.
\end{align*}%\end{subequations}
To establish the second equality, we interchanged the labelling of indices ($p \rightleftarrows q$,\ $k\rightleftarrows \ell$,\ and $m \rightleftarrows n$) 
and we recalled that the partial derivatives commute.%
%We then obtain from the above equation
\end{minipage}
%\flushright
\hspace{-15mm}
\begin{minipage}{0.2\textwidth}
%\begin{figure}[H]
	%\centering
	\unitlength=1.5mm
	\special{em:linewidth 0.4pt}
	\linethickness{0.4pt}
	\begin{picture}(30.00,27.00)%(0,2)
	% wedge
	\put(24.7,10.65){\vector(1,-2){5.00}}
	\put(20.7,5.65){\makebox(0,0)[rc]{\tiny$i$}}
	\put(24.7,10.65){\vector(-1,-2){5.00}}
	\put(28.7,5.65){\makebox(0,0)[lc]{\tiny$j$}}
	% diamond
	\put(19.7,20.65){\vector(1,-2){5.00}}
	\put(20.7,15.65){\makebox(0,0)[rc]{\tiny$k$}}
	\put(29.7,20.65){\vector(-1,-2){5.00}}
	\put(28.7,15.65){\makebox(0,0)[lc]{\tiny$l$}}
	\put(29.7,20.65){\vector(-1,2){5.00}}
	\put(24.7,25.65){\makebox(0,0)[rc]{\tiny$m$}}
	\put(19.7,20.65){\vector(1,2){5.00}}
	\put(25.4,25.65){\makebox(0,0)[lc]{\tiny$n$}}	
	% eyes
	\qbezier(24.7,30.65)(19.7,30.65)(19.7,24.65)
	\put(19.7,24.65){\vector(0,-1){4.00}}
	\put(18.7,25.65){\makebox(0,0)[rc]{\tiny$p$}}
%%%	
	\qbezier(24.7,30.65)(29.7,30.65)(29.7,24.65)
	\put(29.7,24.65){\vector(0,-1){4.00}}
	\put(30.7,25.65){\makebox(0,0)[lc]{\tiny$q$}}
	\end{picture}
%\end{figure}
\end{minipage}
}
Second, the operators encoded using graphs (with a copy of the Poisson bi\/-\/vector~$\cP$ at every internal vertex) can vanish by virtue of the Jacobi identity, see~\eqref{EqJac}, or its differential consequences. 
This mechanism has been illustrated in~\cite{sqs15}; making a part of our present argument 
(see~\cite{f16}), %section~\ref{SecDeform}), 
   %where it is realised in the simplest mode),
it is %will be 
a key to the %(re-)\/
proof of the fact %hypothesis %claim 
that the balanced flow
$\tfrac{\Id}{\Id\veps}(\cP)=\Gamma_1(\cP)+6\,\Gamma_2(\cP)$ does preserve the property of bi\/-\/vectors~$\cP(\veps)$ to be (infinitesimally) Poisson whenever the Cauchy datum~$\cP{\bigr|}_{\veps=0}$ is~such.
\end{rem}
%re\/-\/addressed in full detail in section~\ref{Sec16} on p.~\pageref{Sec16} below.

%Therefore, 
So, each of the two claims (\textit{i}--\textit{ii}) is false if it does not hold for at least one Poisson structure (itself already known to have skew\/-\/symmetric coefficients and turn the left\/-\/hand side of the Jacobi identity into zero for any triple of arguments of the Jacobiator). 
To examine both claims, we %clearly 
need a store of Poisson structures such that the coefficients $\cP^{ij}(\bx)$ are not mapped to zero by the third or second order derivatives in~\eqref{EqFlow1} and~\eqref{EqFlow2}, respectively. For that, a regular generator of Poisson structures with polynomial coefficients of arbitrarily high degree would suffice.

\section{The generators}\label{generators}\noindent%
Let us recall three regular ways to generate %produce 
the Poisson brackets or modify a given one, %which gives 
thus obtaining a new such structure.
%generators of Poisson bi-vectors. 
These generators will be used in section~\ref{Sec12} 
to produce the counterexamples to both claims from~\cite{Ascona96}.

\subsection{The determinant construction}\label{SecDonin}
This generator of Poisson bi-vectors is described in~\cite{Donin}, cf.~\cite{Vinogradov1997} and references therein. The construction goes as follows. Let $x_{1}$,\ $\ldots$,\ $x_{n}$ be the Cartesian coordinates on~$\mathbb{R}^{n\geqslant3}$. Let $\vec{g}=(g_{1}$,\ $\ldots$,\ $g_{n-2})$~be a fixed tuple of smooth functions in these variables. For any $a$,\ $b\in\mathcal{C}^{\infty}(\mathbb{R}^{n})$, put 
\[%\begin{align}
\{a,b\}_{\vec{g}}=\det\bigl(\mathbf{J}(g_{1},\ldots,g_{n-2},a,b)\bigr)
\]%\end{align}
where $\mathbf{J}(\cdot, \ldots, \cdot)$ is the Jacobian matrix. Clearly, the bracket $\{\cdot,\cdot\}_{\vec{g}}$ is bi\/-\/linear and skew\/-\/symmetric. Moreover, it is readily seen to be a derivation in each of its arguments: $\{a,b\cdot c\}_{\vec{g}}=\{a,b\}_{\vec{g}} \cdot c+b \cdot \{a,c\}_{\vec{g}}$. For the validity mechanism of the Jacobi identity for this particular instance of the Nambu bracket we refer to~\cite{Vinogradov1997} again (see also~\cite{Nambu}). 

%To obtain the coefficients $\cP^{ij}(\bx)$ of the respective Poisson bi\/-\/vector~$\cP$, one evaluates the bracket at the coordinate functions: $\cP^{ij}(\bx)=\{x_{i},x_{j}\}_{\vec{g}}{\bigr|}_{\bx}$. 

\begin{example}[see entry~3 in Table~\protect\ref{DonTab} on p.~\protect\pageref{DonTab}]\label{ExDonin}
Fix the functions $g_{1}=x_{2}^{3}x_{3}^{2}x_{4}$ and $g_{2}=x_{1}x_{3}^{4}x_{4}$, and insert them in the determinant generator of Poisson bi\/-\/vectors.
We thus obtain the bi\/-\/vector~$\cP_0$,\label{pExDonin}% 
%the coefficients of which are given in the matrix\marginpar{Order lower indices}
\[
\cP_0^{ij}=\begin{pmatrix}
0&-2\,x_{1}x_{2}^{3}x_{3}^{5}x_{4}&-3\,x_{1}x_{2}^{2}x_{3}^{6}x_{4}&12\,x_{1}x_{2}^{2}x
_{3}^{5}x_{4}^{2}\\ \noalign{\medskip}2\,x_{1}x_{2}^{3}x
_{3}^{5}x_{4}&0&-x_{2}^{3}x_{3}^{6}x_{4}&2\,x_{2}^{3}x_{
	3}^{5}x_{4}^{2}\\ \noalign{\medskip}3\,x_{1}x_{2}^{
	2}x_{3}^{6}x_{4}&x_{2}^{3}x_{3}^{6}x_{4}&0&-3\,x_{2}^{2}
x_{3}^{6}x_{4}^{2}\\ \noalign{\medskip}-12\,x_{1}x_{
	2}^{2}x_{3}^{5}x_{4}^{2}&-2\,x_{2}^{3}x_{3}^{5}x_{4}^{
	2}&3\,x_{2}^{2}x_{3}^{6}x_{4}^{2}&0
\end{pmatrix}.
\]
By construction, the above matrix is skew\/-\/symmetric. The validity of Jacobi identity~\eqref{EqJac} is straightforward: indexed by $i$,\ $j$,\ $k$, all the components $\schouten{\cP,\cP}^{ijk}$ of the tri\/-\/vector vanish.%
\footnote{Indeed, there are four tuples of distinct values of the indices~$i$,\ $j$, and~$k$ up to permutations; we let $1\leqslant i < j < k \leqslant n=4$ so that the check runs over the set of triples $\{(1,2,3)$,\ $(1,2,4)$,\ $(1,3,4)$,\ $(2,3,4)\}$.
%it amounts to the calculation of the tri-vector components $\schouten{\cP,\cP}^{ijk}$ for   . In this case there are $4$, modulo permutations, unique components. Namely:  A quick (computer assisted) computation shows that all these components vanish. 
For example,
%\begin{align*}
$\schouten{\cP,\cP}^{123}%&
=%0+
6x_{1}x_{2}^{5}x_{3}^{11}x_{4}^{2}-6x_{1}x_{2}^{5}x_{3}^{11}x_{4}^{2}
%\\&\phantom{={}}
-6x_{1}x_{2}^{5}x_{3}^{11}x_{4}^{2}+6x_{1}x_{2}^{5}x_{3}^{11}x_{4}^{2}%+0
%\\
%&\qquad%\phantom{={}}
-18x_{1}x_{2}^{5}x_{3}^{11}x_{4}^{2}%+0
+18x_{1}x_{2}^{5}x_{3}^{11}x_{4}^{2}
%\\&\phantom{={}}
+12x_{1}x_{2}^{5}x_{3}^{11}x_{4}^{2}-6x_{1}x_{2}^{5}x_{3}^{11}x_{4}^{2}
-6x_{1}x_{2}^{5}x_{3}^{11}x_{4}^{2}%\\&
=0$.
%\end{align*}
Therefore,
$%\[%\begin{align*}
\schouten{\cP,\cP}=%6
\sum\limits_{1 \leqslant i < j < k \leqslant 4}^{}\schouten{\cP,\cP}^{ijk}(\bx)\,
%\xi_{i}\xi_{j}\xi_{k}
\dd_i\wedge\dd_j\wedge\dd_k =0%.
$.%\]%\end{align*} \marginpar{Edit} %Skew, Poisson, where used below.
} %footnote
This Poisson bi\/-\/vector~$\cP$ is %will be 
used in section~\ref{Sec12} in the list of %our 
counterexamples to the claims under study.
\end{example}

\subsection{Pre\/-\/multiplication in the $3$-\/dimensional case}\label{SecNutku}
%\marginpar{Edit}%Mechanism: Bi-vector <-> one-form.
Let $x,y,z$ be the Cartesian coordinates on the vector space~$\mathbb{R}^3$. For every bi\/-\/vector $\cP=\cP^{ij}\,\partial_{i}\wedge \partial_{j}$, introduce the differential one\/-\/form $\rP=\rP_{1}\,\Id x+\rP_{2}\,\Id y+\rP_{3}\,\Id z$ by setting $\rP\mathrel{{:}{=}}-\cP\inner\Id x\wedge\Id y\wedge\Id z$, so that $\mathrm{P}_{1}=-\cP^{23}$, $\mathrm{P}_{2}=\cP^{13}$, and~$\mathrm{P}_{3}=-\cP^{12}$. %As shown in, 
It is readily seen~\cite{Perelomov} %Nutku} 
that the original Jacobi identity for the bi\/-\/vector~$\cP$ now reads\footnote{%For the one\/-\/form~$\rP$ and its 
The exterior differential~$\Id\rP$ %derivative 
%we have the following expressions:
is equal to
%\begin{subequations}
$%\[%\begin{align*}
%\rP=P_{1}\,\Id x+P_{2}\,\Id y+P_{3}\,\Id z=-\cP^{23}\,\Id x+\cP^{13}\,\Id y-%\cP^{12}\,\Id z,\\
\Id\rP=(\partial_{x} \cP^{13}+\partial_{y} \cP^{23})\,\Id x \wedge \Id y
+(-\partial_{x} \cP^{12}+\partial_{z} \cP^{23})\, \Id x \wedge \Id z
+(-\partial_{y} \cP^{11}-\partial_{z} \cP^{13})\, \Id y \wedge \Id z$.
%\]%\end{align*}
%\end{subequations}
%We can now compute their
The wedge product is %\\[2pt]
%\begin{subequations}
%\mbox{ }\ %
$%\begin{aligned}%\begin{multline*}%
\Id\rP\wedge\rP%&
= 
\bigl(\partial_{x}\cP^{31}\,\cP^{12}+\partial_{y}\cP^{23}\,\rP^{21}
+\partial_{x}\cP^{12}\,\cP^{13}+\partial_{z}\cP^{23}\,\cP^{31}
+\partial_{y}\cP^{12}\,\cP^{23}+\partial_{z}\cP^{31}\,\cP^{32}\bigr)\,\Id x \wedge\Id y \wedge \Id z%\\ 
%&
= (-\schouten{\cP,\cP}\inner \Id x\wedge\Id y\wedge \Id z)\,\Id x\wedge\Id y\wedge \Id z$.
%\end{aligned}$%\end{multline*}%
%\end{subequations}
} %footnote
$\Id\rP\wedge\rP=0$ for the respective one\/-\/form~$\rP$. %Therefore, 
But let us note that the pre\/-\/multiplication $\rP\mapsto f\cdot\rP$ of the form~$\rP$ by a smooth function~$f$ preserves this reading of the Jacobi identity:
$%\[%\begin{align}
\Id(f\rP)\wedge(f\rP) = f \cdot \bigl[\Id f \wedge\rP \wedge\rP + f\cdot\Id\rP \wedge\rP \bigr] = f^{2}\cdot\Id\rP \wedge\rP=0$.
%\]%\end{align}
This shows that %whence 
the bi\/-\/vector~$f\cP$ which the form $f\rP$ yields on~$\mathbb{R}^{3}$ is also Poisson.

This pre\/-\/multiplication trick provides the examples of Poisson structures of arbitrarily high polynomial degree coefficients (in a manifestly non\/-\/symplectic three\/-\/dimensional set\/-\/up).\footnote{\label{flownonzero}%
In dimension three, this pre\/-\/multiplication procedure also provides the examples of Poisson bi\/-\/vectors at which the second flow~\eqref{EqFlow2} does not vanish identically.}
%, and are therefore counterexamples. }

\subsection{The Vanhaecke construction}\label{SecVanhaecke}
In~\cite{Van}, Vanhaecke created another construction of high polynomial degree Poisson bi\/-\/vectors. Let $u$~be a monic degree~$d$ polynomial in~$\lambda$ and $v$~be a polynomial of degree $d-1$ in~$\lambda$:
\[%\begin{align*}
u(\lambda)%&
=\lambda^d+u_{1}\lambda^{d-1}+\ldots+u_{d-1}\lambda+u_{d},\qquad%\\
v(\lambda)%&
=%\phantom{\Lambda^d+{}}
v_1\lambda^{d-1}+\ldots+v_{d-1}\lambda+v_{d}.
\]%\end{align*}
Consider the space $\Bbbk^{2d}$ (e.g., set $\Bbbk\mathrel{{:}{=}}\mathbb{R}$) with Cartesian coordinates $u_{1}$,\ $\ldots$,\ $u_{n}$,\ $v_{1}$,\ $\ldots$,\ $v_{d}$. 
To define the Poisson bracket, fix a bivariate polynomial~$\phi(\cdot,\cdot)$ %over variables $x$ and $y$ which will be evaluated in $u_{i}$'s and $v_{j}$
and for all $1\leqslant i,j\leqslant d$ set
%\begin{subequations}
\begin{equation}\label{VanHaeckePoisson}
%\begin{align}
\{u_{i},u_{j}\}%&
=\{v_{i},v_{j}\}=0,\qquad% \\
\{u_{i},v_{j}\}%&
={}\text{coeff.\ of $\lambda^j$ in}\ %\quad
%\left(
\phi\bigl(\lambda,v(\lambda)\bigr)\cdot\left[\frac{u(\lambda)}{\lambda^{d-i+1}}\right]_{+} \mod u(\lambda)%\right)
, %_{\lambda^{j}} 
%\end{align}
\end{equation}
%\end{subequations}
where we denote by $[\ldots]_{+}$ the argument's polynomial part
% of the argument, i.e. discarding the rational part, 
and where the remainder mo\-du\-lo the degree~$d$ polynomial~$u(\lambda)$ is obtained %by 
using the Euclidean division algorithm.
%mod implies performing  and $[\ldots]_{\lambda^{j}}$ means taking the coefficient of the term $\lambda^{j}$. 
%After the Euclidean division, 
     %%We are guaranteed by the Euclidean division that 
%the highest power of~$\lambda$ in the remainder is not exceeding~$d-1$.%\marginpar{Edit}%Next phrase is a tautology.
%meaning that the obtained coefficients belong to the monomials $\lambda^{d-1}$ to $\lambda^{0}$. 

Let us emphasise that these Poisson bi\/-\/vector are defined on the even\/-\/dimensional spaces. Indeed, the coefficients of Poisson bracket~\eqref{VanHaeckePoisson} are arranged in the block matrix $\left(\begin{smallmatrix}\phantom{+}0 & U\\ -U&0\end{smallmatrix}\right)$, where the components of the matrix~$U$ are~$U^{ij}=\{u_{i},v_{j}\}$.
%We can now write the bi-vector for this Poisson bracket using the above coefficients in order to obtain a  of the form: ,  One should notice that this 

\subsection{The Hamiltonian differential operators on jet spaces}\label{SecHD}
The variational Poisson brackets $\{\cdot,\cdot\}_{\boldsymbol{\cP}}$ for functionals of sections of affine %fibre 
bundles ge\-ne\-ra\-li\-se the notion of Poisson brackets $\{\cdot,\cdot\}_{\cP}$ for functions on finite\/-\/dimensional Poisson manifolds $({N}^{n},\{\cdot,\cdot\}_{\cP})$. Namely, let us consider the space~${J}^{\infty}(\pi)$ of infinite jets of sections for a given bundle~$\pi$ over a manifold~${M}^{n}$ of positive dimension~$m$. The variational Poisson brackets $\{\cdot, \cdot\}_{\boldsymbol{\cP}}$ on~${J}^{\infty}(\pi)$ are then specified by using the Hamiltonian differential operators (which we shall denote by~$A$ and the order of which is typically positive).\footnote{In fact, the Poisson geometry of finite\/-\/dimensional affine manifolds~$({N}^{n},\{\cdot,\cdot\}_{\cP})$ is a zero differential order sub\/-\/theory in the variational Poisson geometry of infinite jet spaces~${J}^{\infty}(\pi)$. Indeed, let the fibres in the bundle~$\pi$ be~${N}^{n}$ and proclaim that only \emph{constant} sections are allowed.}
The formalism of variational Poisson bi\/-\/vectors $\boldsymbol{\cP}=\frac{1}{2}\langle \boldsymbol{\xi}\cdot \vec{A}(\boldsymbol{\xi}) \rangle$ and the variational Schouten bracket $\schouten{\cdot,\cdot}$ is standard (see~\cite{Olver,cycle14}% and~\cite{f16}% 
%section~\ref{SecDeform} below
). The geometry of iterated variations is revealed in~\cite{gvbv}; 
the correspondence between the Kontsevich graphs and local variational polydifferential operators is explained in~\cite{dq15}.

\begin{example}
%For an 
To inspect whether either of the two claims (which we quote from~\cite{Ascona96} on the title page) would hold in the variational set\/-\/up, it is enough to consider a Hamiltionian differential operator with 
(differential-)\/polynomial coefficients of degree~$\geqslant 3$. 
Let us %conveniently 
take the Hamiltonian operator\footnote{More examples of variational Poisson structures, which are relevant for our present purpose, can be found in~\cite{Olver1997} or, e.g., in~\cite{Vodova} (see also the references contained therein).}
$%\[%\begin{align}
A = u^{2} \circ %\frac
{\Id}/{\Id x} \circ u^{2}
$ %\]%\end{align}
for the Harry Dym equation (see \cite{Kent}); here $u$~is the fibre coordinate in the trivial bundle~$\pi\colon \mathbb{R}\times \mathbb{R}\rightarrow \mathbb{R}$ and $x$~is the base variable. This operator is obviously skew\/-\/adjoint, whence the variational Poisson bracket $\{\cdot,\cdot\}_{\boldsymbol{\cP}}$ is skew\/-\/symmetric. The Jacobi identity for $\{\cdot,\cdot\}_{\boldsymbol{\cP}}$ is also easy to check: the variational master\/-\/equation $\schouten{\boldsymbol{\cP},\boldsymbol{\cP}} \cong 0$ does hold for the variational bi\/-\/vector $\boldsymbol{\cP}=\frac{1}{2}\langle 
%\boldsymbol
{\xi}\cdot \vec{A}(%\boldsymbol
{\xi}) \rangle$.
\end{example}

\section{The counterexamples}\label{Sec12}\noindent%
We now examine the properties of both tetrahedral flows~\eqref{EqFlows} whenever each of them is evaluated at a given Poisson bi\/-\/vector. (Examples of such bi\/-\/vectors are produced by using the techniques from section~\ref{generators}.) To motivate the composition of Tables~\ref{DonNutTab}\/--\/\ref{HDTab} and clarify the meaning of their content, let us consider %first 
an example: namely, we first take the Poisson bi\/-\/vector which was obtained in section~\ref{SecDonin} (see p.~\pageref{pExDonin}).
%We can now consider the evaluation of the flows, corresponding to the tetrahedron, on Poisson bi-vectors obtained from the generators in section $2$. We first  explicitly write out the results of the example described in section \ref{generators}. 

\begin{example}[continued]\label{ExContinued}%\marginpar{Order lower indices}
Rewriting the Poisson bi\/-\/vector $\cP_0\in\Gamma\bigl(\bigwedge^2 TN^4\bigr)$ in terms of the parity\/-\/odd variables~$\boldsymbol{\xi}$, we obtain that under the isomorphism $\Gamma\bigl(\bigwedge^\bullet TN^n\bigr) \simeq C^\infty(\Pi T^*N^n)$ the bi\/-\/vector $\cP_0^{ij}(\bx)\,\dd_i\wedge\dd_j$ becomes~$\tfrac{1}{2}\cP_0^{ij}(\bx)\,\xi_i\xi_j$;
we have that $\cP_{0}=$
%The Poisson bi-vector from section~\ref{generators} can be written down as a sum of components, rather than a matrix, yielding the expression
\[%\begin{align*}
%&=
-2\,x_{1}x_{2}^{3}x_{3}^{5}x_{4}\xi_{1}\xi_{2}
-3\,x_{1}x_{2}^{2}x_{3}^{6}x_{4}\xi_{1}\xi_{3}
+12\,x_{1}x_{2}^{2}x_{3}^{5}x_{4}^{2}\xi_{1}\xi_{4}%\\
%&\qquad%\phantom{{}=}
-x_{2}^{3}x_{3}^{6}x_{4}\xi_{2}\xi_{3}
+2\,x_{2}^{3}x_{3}^{5}x_{4}^{2}\xi_{2}\xi_{4}
-3\,x_{2}^{2}x_{3}^{6}x_{4}^{2}\xi_{3}\xi_{4}.
\]%\end{align*}
Now, we calculate the right\/-\/hand sides $\cP_{1}\mathrel{{:}{=}}\Gamma_1(\cP_0)$ and $\cP_{2}\mathrel{{:}{=}}\Gamma_2(\cP_0)$ of tetrahedral flows~\eqref{EqFlows}.
%flows corresponding to each of the tetrahedra, thereby obtaining bi-vectors $\cP_{1}$ and $\cP_{2}$. 
The coefficient matrix of the bi\/-\/vector~$\cP_{1}$ is 
%obtained from evaluating the flow corresponding to $\Gamma_{1}$ into the above Poisson bi-vector:
\begin{align*}
\cP_{1}^{ij}= \begin{pmatrix}
0&-24480\,x_{1}x_{2}^{9}x_{3}^{20}x_{4}^{4}&-51840\,x_{1
}x_{2}^{8}x_{3}^{21}x_{4}^{4}&12960\,x_{1}x_{2}^{8}x_{3}^{20}x_{4}^{5}
\\ \noalign{\medskip}24480\,x_
{1}x_{2}^{9}x_{3}^{20}x_{4}^{4}&0&-15480\,x_{2}^{9}x_{3}^{21}x_{4}^{4}&2448\,x_{2}^{9}x_{3}
^{20}x_{4}^{5}\\ \noalign{\medskip}51840\,x_{1}x_{2}^{8}x_{3}
^{21}x_{4}^{4}&15480\,x_{2}^{9}x_{3}^{21}x_{4}^
{4}&0&-18144\,x_{2}^{8}x_{3}^{21}x_{4}^{5}
\\ \noalign{\medskip}-12960\,x
_{1}x_{2}^{8}x_{3}^{20}x_{4}^{5}&-2448\,x_{2}^{9}x_{3}^{20}x_{4}^{5}&18144\,x_{2}^{8}x_{3}
^{21}x_{4}^{5}&0
\end{pmatrix}.
\end{align*}
In a similar way, % Similarly, 
the polydifferential operator~$\Gamma_2$ (encoded by the second tetrahedral graph in Fig.~\ref{FigTetra}) yields the matrix
%bi-vector $\cP_{2}$ is obtained from evaluating the flow corresponding to $\Gamma_{2}$:
\begin{align*}
\cP_{2}^{ij}= \begin{pmatrix}
16920x_{1}^{2}x_{2}^{8}x_{3}^{20}x_{4}
^{4}&-12060\,x_
{1}x_{2}^{9}x_{3}^{20}x_{4}^{4}&-16380\,x_{1}x_{2}^{8}x_{3}^{21}x_{4}^{4}&42840\,
x_{1}x_{2}^{8}x_{3}^{20}x_{4}^{5}\\ \noalign{\medskip}
2700\,x_{1}x_{2}^{9}x_{3}^{20}x_{4}^{4}&-7200\,x_{2}
^{10}x_{3}^{20}x_{4}^{4}&4680\,x_
{2}^{9}x_{3}^{21}x_{4}^{4}&-252\,x_{2}^{9}x_{3}^{20}x_{4}^{5}
\\ \noalign{\medskip}-13140\,x
_{1}x_{2}^{8}x_{3}^{21}x_{4}^{4}&5040\,x_{2}^{9}x_{3}^{21}x_{4}^{4}&-12060\,x_{2}
^{8}x_{3}^{22}x_{4}^{4}&13716\,x_
{2}^{8}x_{3}^{21}x_{4}^{5}\\ \noalign{\medskip}-80280\,x_{1}x_
{2}^{8}x_{3}^{20}x_{4}^{5}&-18036\,x_{2}^{9}x_{3}^{20}x_{4}^{5}&
21708\,x_{2}^{8}x_{3}^{21}x_{4}^{5}&-58104\,x_{2}^{8}x_{3}^{20}
x_{4}^{6}
\end{pmatrix}.
\end{align*}
Notice that this coefficient matrix %corresponding to the bi-vector $\cP_{2}$
is not yet antisymmetric, 
%yet the Poisson bi-vector should by definition always be anti-symmetric.
but its \emph{symmetric} counterpart is skipped out %anyway 
in the construction of the bi\/-\/vector~$\cP_2$ and its transcription by using the anticommuting variables~$\boldsymbol{\xi}$.
Therefore, we antisymmetrise the above matrix at once, the output to be used in what follows. %for further computations. 
We obtain that the bi\/-\/vector %~$\cP_{2}$ 
is
%The anti-symmetrised matrix for the  expressed in components is:
\begin{align*}
\cP_{2} &= 
-7380x_{1}x_{2}^{9}x_{3}^{20}x_{4}^{4}  \xi_{1}\xi_{2}-1620x_{1}x_{2}^{8}x_{3}^{21}x_{4}^{4}  \xi_{1}\xi_{3}+61560x_{1}x_{2}^{8}x_{3}^{20}x_{4}^{5}  \xi_{1}\xi_{4} 
\\&\phantom{{}=}
-180x_{2}^{9}x_{3}^{21}x_{4}^{4}  \xi_{2}\xi_{3}+8892x_{2}^{9}x_{3}^{20}x_{4}^{5}  \xi_{2}\xi_{4}-3996x_{2}^{8}x_{3}^{21}x_{4}^{5}  \xi_{3}\xi_{4}.
\end{align*}
We now see that for the Poisson bi\/-\/vector~$\cP_0$ from Example~\ref{ExDonin} on p.~\pageref{pExDonin}, \textit{\textbf{the bi\/-\/vector~$\cP_2$ does not vanish}},
%One should now notice that the bi-vector $\cP_{2}$ does not vanish, 
thereby disavowing %contradicting 
the second claim from~\cite{Ascona96}.

To check the compatibility of the original Poisson bi\/-\/vector~$\cP_{0}$ with the newly obtained bi-vector $\cP_{1}$, we calculate %evaluate 
their Schouten bracket:%\marginpar{Edit}%Divide by 2.
\begin{align*}
%\tfrac{1}{3!}
\schouten{\cP_{0},\cP_{1}}&= 
46008\,x_{1}x_{2}^{11}x_{3}^{26}x_{4}^{5}\xi_{1}\xi_{2}\xi_{3}
+852768\,x_{1}x_{2}^{11}x_{3}^{25}x_{4}^{6}\xi_{1}\xi_{2}\xi_{4}\\
&\qquad%\phantom{{}=}
+1246752\,x_{1}x_{2}^{10}x_{3}^{26}x_{4}^{6}\xi_{1}\xi_{3}\xi_{4}
+340200\,x_{2}^{11}x_{3}^{26}x_{4}^{6}\xi_{2}\xi_{3}\xi_{4}\neq0.
%[P_{0},P_{1}] =648(71\,x_{2}^{11}x_{3}^{26}x_{4}^{5}x_{1}b_{1}b_{2}b
%_{3}+1316\,x_{1}x_{2}^{11}x_{3}^{25}x_{4}^{6}b_{1}
%b_{2}b_{4}+\\1924\,x_{1}x_{2}^{10}x_{3}^{26}x_{4}^{
%6}b_{1}b_{3}b_{4}+525\,x_{2}^{11}x_{3}^{26}x_{4}^
%{6}b_{2}b_{3}b_{4})
\end{align*}
The above expression is not identically zero. Therefore, %we conclude that 
\textit{\textbf{the leading %first flow, as a first 
term~$\cP_1$ in the deformation %in the expression 
$\cP_0\mapsto\cP(\varepsilon)=\cP_0+\varepsilon\cP_1+\bar{o}(\varepsilon)$
destroys the property of bi\/-\/vector~$\cP(\varepsilon)$ to be Poisson}} at~$\varepsilon\neq0$ on all of~$\BBR^4$.
%for all~$\bx\in\BBR^4$.
%, does not make the deformation  . 

%Just for the sake of curiosity, let us check whether the 
The same compatibility test for~$\cP_0$ and its second flow~\eqref{EqFlow2} yields that %%\marginpar{Edit}%Divide by 2.
%of the flow corresponding to $\Gamma_{2}$ holds, by evaluating the Schouten bracket
\begin{align*}
%\tfrac{1}{3!}
\schouten{\cP_{0},\cP_{2}}&=
-{7668}\,x_{1}x_{2}^{11}x_{3}^{26}x_{4}^{5}\xi_{1}\xi_{2}\xi_{3}
-{142128}\,x_{1}x_{2}^{11}x_{3}^{25}x_{4}^{6}\xi_{1}\xi_{2}\xi_{4}\\
&\qquad%\phantom{{}=}
-{207792}\,x_{1}x_{2}^{10}x_{3}^{26}x_{4}^{6}\xi_{1}\xi_{3}\xi_{4}
-{56700}\,x_{2}^{11}x_{3}^{26}x_{4}^{6}\xi_{2}\xi_{3}\xi_{4}.
%[P_{0},P_{2}] =108(-71\,x_{2}^{11}x_{3}^{26}x_{4}^{5}x_{1}b_{1}b_{2}b
%_{3}-1316\,x_{1}x_{2}^{11}x_{3}^{25}x_{4}^{6}b_{1}
%b_{2}b_{4}-\\1316\,x_{1}x_{2}^{10}x_{3}^{26}x_{4}^{6
%}b_{1}b_{3}b_{4}-525\,x_{2}^{11}x_{3}^{26}x_{4}^{6
%}b_{2}b_{3}b_{4})
\end{align*}
Again, this expression %for the evaluated Schouten bracket 
does not vanish identically on all %at points~$\bx$ 
of the Poisson manifold $\bigl(\BBR^4$,\ $\{\cdot,\cdot\}_{\cP_0}\bigr)$.
%%%
%Therefore we conclude that 
%neither the first nor the second flow preserve the Poisson property of $\cP$.
%%%
We conclude that %, so that 
neither %$\cP_{1}$ nor $\cP_{2}$ 
of two flows~\eqref{EqFlows}
preserve the property of bi\/-\/vector~$\cP(\varepsilon)$ to stay (infinitesimally) Poisson at~$\varepsilon\neq0$ for this example of Poisson bi\/-\/vector.%when chosen as the deformation leading term $\cQ$.
\footnote{Let us also inspect %We can now verify 
whether the Jacobi identity holds for any of the %obtained 
bi\/-\/vectors~$\cP_{1}$ and~$\cP_{2}$.
For~$\cP_{1}$ we have that the left\/-\/hand side of the Jacobi identity is equal to%\marginpar{Edit}%Divide by 2.
\[%\begin{align*}
\schouten{\cP_{1},\cP_{1}} =-2963589120\cdot %3!\cdot
\bigl(x_{1}x_{2}^{17}x_{3}^{41}x_{4}^{8}\xi_{1}\xi_
{2}\xi_{3}+5\,x_
{1}x_{2}^{17}x_{3}^{40}x_{4}^{9}\xi_{1}\xi_{2}\xi_{4}-2\,x_{1}x_{
	2}^{16}x_{3}^{41}x_{4}^{9}\xi_{1}\xi_{3}\xi_{4}\bigr),
\]%\end{align*}
which does not vanish. 
%(Therefore the Jacobi identity is not satisfied for~$\cP_1$.) 
For~$\cP_{2}$ the left\/-\/hand side of the Jacobi identity equals%\marginpar{Edit}%Divide by 2.
\[%\begin{align*}
\schouten{\cP_{2},\cP_{2}} =-262517760\cdot %3!\cdot
\bigl(x_{1}x_{2}^{17}x_{3}^{41}x_{4}^{8}\xi_{1}\xi_
{2}\xi_{3}+5\,x_{1
}x_{2}^{17}x_{3}^{40}x_{4}^{9}\xi_{1}\xi_{2}\xi_{4}-2\,x_{1}x_{2}
^{16}x_{3}^{41}x_{4}^{9}\xi_{1}\xi_{3}\xi_{4}\bigr).
\]%\end{align*}
This expression also does not vanish, so that neither~$\cP_{1}$ nor~$\cP_{2}$ are Poisson bi\/-\/vectors.}
\end{example}

\begin{rem}%We do notice that
In %for this 
the above example, the %evaluated 
Schouten brackets $\schouten{\cP_{0},\cP_{1}}$ and $\schouten{\cP_{0},\cP_{2}}$ are determined by %consist of 
the same polynomials in the variables~$\bx$ and~$\boldsymbol{\xi}$:
%  Indeed, a quick check shows 
we see that $\schouten{\cP_{0},\cP_{1}}=-6\cdot \schouten{\cP_{0},\cP_{2}}$. 
This implies that for this example of Poisson bi\/-\/vector~$\cP_0$,
the leading term $\cQ\mathrel{{:}{=}}\cP_{1}+6\cP_{2}$ does (infinitesimally) preserve the property of $\cP(\varepsilon)$ to be Poisson
in the course of deformation $\cP_0\mapsto\cP_0+\varepsilon\cQ+\bar{o}(\varepsilon)$.% , by choosing .  That way 
%Indeed, we have that
%\[\schouten{\cP_{0},\cQ}=\schouten{\cP_{0},\cP_{1}+6\cP_{2}}=\schouten{\cP_{0},\cP_{1}}+6\,\schouten{\cP_{0},\cP_{2}}=0\]
%due to the linearity of the Schouten bracket. 

Moreover, it is readily seen that the ratio $1:6$ is the \emph{only} way to balance the two flows,~\eqref{EqFlow1} vs~\eqref{EqFlow2}, such that their nontrivial linear combination~$\cQ$ is compatible with the %given 
Poisson bi\/-\/vector~$\cP_0$ from Example~\ref{ExDonin}.\footnote{The balance $1:\tfrac{4}{3}$ was considered %advocated 
in~\cite[\S5.2]{Merkulov} for the linear combination of flows~\eqref{EqFlow1} and~\eqref{EqFlow2}, respectively%
%; our present argument and the counterexamples which follow withdraw that claim
.}%
\end{rem}

\begin{rem}%Note that 
In Example~\ref{ExContinued}
the linear combination~$\cQ=\cP_{1}+6\cP_{2}\neq0$ of two flows~\eqref{EqFlows} 
is not identically equal to zero.
(For other examples this may happen incidentally.)
The leading term~$\cQ$ in the infinitesimal deformation $\cP_0\mapsto\cP_0+\veps\cQ+\bar{o}(\veps)$ is %non
trivial in the Poisson cohomology with respect to~$\boldsymbol{\dd}_{\cP_0}$, i.\,e.\ %that is, 
$\cQ=%\neq
\schouten{\cP_0,\cX}$ for some %any 
vector~$\cX$ on the four\/-\/dimensional %vector 
space.\footnote{In all the two\/-\/dimensional Poisson geometries, 
%It is shown in~\cite{Ascona96} that in the 2-dimensional case the first flow 
the first flow~$\cP_{1}$ is always cohomologically trivial, i.e.\ it is of the form $\cP_{1}=\schouten{\cP_0,\cX}$ for some one\/-\/vector~$\cX$, see~\cite{Ascona96}.}
Hence this %bi\/-\/vector
$\cQ$~is trivially compatible with
the %initial 
Poisson bi\/-\/vector~$\cP_{0}$: namely, $\lshad\cP_0,\cQ\rshad\equiv0$, 
see p.~\pageref{pConclusion} below.
\end{rem}
%can be %need not be 
%compatible with in a nontrivial manner, that is, the  %does not have to be
%, although 
%In the above example one can check that $\cQ=\cP_{1}+6\cP_{2} \neq 0$. 
%We expect %also note 

\smallskip%
In the three tables below we summarise the results about the 
flows~$\cP_{1}$ and~$\cP_{2}$, which we evaluate at the examples of Poisson bi\/-\/vectors~$\cP_0$. Special attention is paid to the leading deformation term~$\cQ=\cP_1+6\cP_2$ in each case: we inspect whether this bi\/-\/vector incidentally vanishes and whether it is (indeed, always) compatible with the original Poisson structure~$\cP_0$.
%and $\cQ$ in the Poisson bi-vectors obtained from the methods described in section~\ref{generators}, are summarised in the tables below.

\begin{table}[htb]
\centering
\caption{The Poisson bi\/-\/vectors~$\cP_{0}$ are generated using the determinant method from section~\ref{SecDonin} (the dimension is equal to~$3$, so we specify the fixed argument~$g_{1}$); that generator is combined with the pre\/-\/multiplication~$(f\cdot)$ as explained in section~\ref{SecNutku}.}%
\label{DonNutTab}
\vspace{1mm}
\begin{tabular}{|l|c|l|c|c|c|c|c|}
\hline
\textnumero  & dim & Argument \& pre-factor &  $\schouten{\cP_{0},\cP_{1}}$ & 
$\cP_{2}\stackrel{\text{?}}{=}0$\rule{0pt}{5mm} & $\schouten{\cP_{0},\cP_{2}}$ & 
$\cQ\stackrel{\text{?}}{=}0$ & $\schouten{\cP_{0},\cQ}$
\\
 &  &  &  \quad${}= 0$\,? &  &\quad${}=0$\,? &  & \quad${}= 0$\,? \\
\hline
1.\rule{0pt}{6mm} 
& 3 & $[x_{1}^{5}x_{2}^{3}x_{3}^{4}+x_{1}^{2}x_{3}^{5}+x_{1}x_{2}^{5}x_{3}]$ &  \xmark & \xmark & \xmark  & \xmark & \vmark \\
 &  & ${x_{1}^3+x_{2}^2}$   &  &  & & &\\
2. & 3 & $[x_{1}x_{2}+x_{1}x_{3}+x_{2}x_{3}]$ &  \xmark & \xmark & \xmark & \xmark & \vmark \\
 &  & ${x_{1}^2+x_{2}}$ &    &  & & & \\[1mm]
\hline
%3 & 3 & $[x_{1}x_{2}+x_{3}x_{1}]$ & \xmark & \xmark & \xmark & \xmark \\
% &  & $l_{y_{1}^3+y_{2}}$ (Recheck this!)&  &  &  & \\
\end{tabular}
\end{table}
For both %the above 
examples in Table~\ref{DonNutTab}
we have that neither does~$\cP_{1}$ preserve the property of~$\cP_{0}+\varepsilon \cP_{1}+\bar{o}(\varepsilon)$ to be (infinitesimally) Poisson nor does~$\cP_{2}$ vanish identically~--- which is in contrast with both the claims from~\cite{Ascona96}.

\begin{table}[htb]
\centering
\caption{In dimensions higher than $3$, we %now 
generate the Poisson bi\/-\/vec\-tors~$\cP_{0}$ by using the determinant method from section~\ref{SecDonin}: the auxi\-li\-ary arguments~$g_{1}$,\ $\ldots$,\ $g_{n-2}$ are specified.}\label{DonTab}
\vspace{1mm}
\begin{tabular}{|l|c|l|c|c|c|c|c|}
	\hline
	\textnumero  & dim & Arguments & $\schouten{\cP_{0},\cP_{1}}$ & 
	$\cP_{2}\stackrel{\text{?}}{=}0$\rule{0pt}{5mm} & $\schouten{\cP_{0},\cP_{2}}$ & 
	$\cQ\stackrel{\text{?}}{=}0$ & $\schouten{\cP_{0},\cQ}$ \\
	&  &  &  \quad${}=0$\,? &  &\quad${}=0$\,? &  & \quad${}=0$\,?\\
	\hline
	3.\rule{0pt}{6mm}  
	& 4 & $[x_{2}^{3}x_{3}^{2}x_{4},x_{1}x_{3}^{4}x_{4}]$ &  \xmark & \xmark & \xmark & \xmark & \vmark\\
	4.  &4 & $[x_{1}^{2}x_{2}^{3}x_{3}^{4}x_{4}^{5},x_{1}x_{2}x_
	{3}x_{4}]$ & \xmark & \xmark & \xmark & \vmark & \vmark\\
	5.  & 4 & $[x_{2}^{2}x_{3}^{2}x_{4}^{2},x_{1}^{2}x_{3}^{2}x_{4}^{2}]$  & \xmark & \xmark & \xmark & \vmark & \vmark\\
	6.  & 5 & $[x_{2}^{3}x_{3}^{2}x_{4},x_{1}x_{3}^{4}x_{4},x_{3}^{3}x_{4}^{2}x_{5}^{4}]$  & \xmark & \xmark & \xmark & \xmark & \vmark\\[1mm]
	\hline
\end{tabular}
\end{table}
In Table~\ref{DonTab} we again have that 
neither is the property to be (infinitesimally) Poisson 
preserved for $\cP_{0}+\varepsilon \cP_{1}+\bar{o}(\varepsilon)$ 
nor is the bi\/-\/vector~$\cP_{2}$ vanishing identically. 
%The flows~$\cP_{1}$ and~$\cP_{2}$ %are 
%can occasionally be Poisson: e.g., see example~5. 
%The combined flow $\cQ$ may also occasionally be Poisson, e.g., example 2.

\begin{table}[htb]
\centering
\caption{The results for the Vanhaecke method from section~\ref{SecVanhaecke}:
we here specify the bivariate polynomials~$\phi$.%
%using as argument the function $\phi(x,y)$.
}\label{VanTab}
\vspace{1mm}
\begin{tabular}{|l|c|l|c|c|c|c|c|}
\hline
\textnumero  & dim & $\phi(x,y)$ &  $\schouten{\cP_{0},\cP_{1}}\stackrel{\text{?}}{=}0$ & 
$\cP_{2}\stackrel{\text{?}}{=}0$\rule{0pt}{5mm} & $\schouten{\cP_{0},\cP_{2}}\stackrel{\text{?}}{=}0$  & $\cQ\stackrel{\text{?}}{=}0$ & $\schouten{\cP_{0},\cQ}\stackrel{\text{?}}{=}0$\\
\hline
7.\rule{0pt}{6mm}
  & 4 & $[x^{2}y^{2}]$ & \xmark & \xmark & \xmark & \xmark & \vmark\\
8.  & 4 & $[x^{2}y]$  & \xmark & \xmark & \xmark & \xmark & \vmark\\
9.  & 4 & $[x^{3}y^{2}]$  & \xmark & \xmark & \xmark & \xmark & \vmark\\
10.  & 4 & $[x^{3}y^{3}]$  & \xmark & \xmark & \xmark & \xmark & \vmark\\
11.  & 6 & $[x^{2}y^{2}]$  & \xmark & \xmark & \xmark & \xmark & \vmark\\[1mm]
\hline
\end{tabular}
\end{table}
The entries in Table~\ref{VanTab} report on the use of the generator from section~\ref{SecVanhaecke}:
experimentally established, the properties of %For 
these Poisson bi\/-\/vectors do not match %we also have that they contradict with 
both the claims from~\cite{Ascona96}.

\begin{table}[htb]
\centering
\caption{The results for the infinite\/-\/dimensional case.}\label{HDTab}
\vspace{1mm}
\begin{tabular}{|l|c|l|c|c|c|c|%c|
}
\hline
\textnumero  & dim & Operator &  $\schouten{\boldsymbol{\cP}_{0},\boldsymbol{\cP}_{1}}\stackrel{\text{?}}{\cong}0$ & $\boldsymbol{\cP}_{2}\stackrel{\text{?}}{\cong}0$\rule{0pt}{5mm}%& $\schouten{\boldsymbol{\cP}_{0},\boldsymbol{\cP}_{2}}\stackrel{\text{?}}{\cong}0$
\\
\hline
12.\rule{0pt}{6mm}  
& $\infty$ & $u^{2}\circ %\tfrac
{\Id}/{\Id x} \circ u^{2}$ &  \xmark & \vmark% 
% (!): AVK 31/03/2016.    & \xmark
\\[1mm]
\hline
\end{tabular}
\end{table}
%
%HDym flows:
The variational bi\/-\/vector~$\boldsymbol{\cP}_1=\tfrac{1}{2}\langle\xi\cdot\vec{A}_1(\xi)\rangle$, which we construct from the variational Poisson bi\/-\/vector $\boldsymbol{\cP}_0=\tfrac{1}{2}\langle\xi\cdot u^2\,%\tfrac
{\vec{\Id}}/{\Id x}(u^2\,\xi)\rangle$ by using the geometric technique from~\cite{gvbv} (see also~\cite{dq15}), is determined by the (skew\/-\/adjoint part of the) %\marginpar{Edit} %skew\/-\/adjoint
first\/-\/order differential operator $A_1=192\,\bigl(9u^{8}u_{x}u_{xx}-u^{9}u_{xxx}\bigr)\,%\tfrac
{\Id}/{\Id x}$ in total derivatives. Again (see Table~\ref{HDTab}), the two variational bi\/-\/vectors are \emph{not} compatible: we check that $\schouten{\boldsymbol{\cP}_0,\boldsymbol{\cP}_1} \ncong0$ under the %for their 
variational Schouten bracket.
%%%
Remarkably, %Likewise, 
the variational bi\/-\/vector~$\boldsymbol{\cP}_{2}$ is specified by the 
%self\/-\/adjoint
second\/-\/order total dif\-fe\-ren\-tial operator 
%$A_2={\ }\neq0$, %\marginpar{Edit}
whose skew\/-\/adjoint component vanishes, 
%=-(2592(uu_{xx}+u_{x}^{2}))u^{6}u_{x}^{2}+5184u^{6}u_{x}^{4}+0D_{x}+(288u^{8}(uu_{xx}+u_{x}^{2})-576u^{8}u_{x}^{2})D_{x}^{2}
whence the respective variational bi\/-\/vector \emph{is} %not 
equal to zero (modulo exact terms 
%up to any integration by parts 
within its horizontal cohomology class~\cite{Olver}).

%\section{The $1:6$ balance of flows}\label{Sec16}
\section*{Conclusion}\label{pConclusion}
The linear combination $\cQ=\cP_{1}+6\cP_{2}$ of the Kontsevich tetrahedral flows preserves the space of Poisson bi\/-\/vectors~$\cP_{0}$ under the infinitesimal deformations $\cP_{0}\mapsto\cP_{0}+ \varepsilon \cQ + \bar{o}(\varepsilon)$. This is manifestly true for all the examples of Poisson 
bi\/-\/vectors on finite\/-\/dimensional (vector or affine) spaces~$\mathbb{R}^{n}$ which we have considered so far. 
We conjectured that the leading deformation term~$\cQ=\cQ(\cP_{0})$ always has this property, that is, 
\textit{\textbf{the bi\/-\/vector~$\cQ$ marks a $\boldsymbol{\partial}_{\cP_0}$-\/cohomology class for every Poisson bi\/-\/vector~$\cP_{0}$}} on a finite\/-\/dimensional affine manifold. (Recall that such class can be $\boldsymbol{\partial}_{\cP_0}$-\/trivial; moreover, %it 
the bi\/-\/vector~$\cQ$ can vanish identically~--- yet the above examples confirm the existence of Poisson geometries where neither of the two options is realised.)

Let us conclude that every claim of an object's vanishing by virtue of the skew\/-\/symmetry and Jacobi identity for a given Poisson bi\/-\/vector, which that object depends on by construction, must be accompanied with an explicit description of that factorisation mechanism (e.g., see~\cite{sqs15}) or at least, with a proof of that mechanism's existence. Apart from the trivial case (here, $\cQ=0$ so that $\schouten{\cP_{0},\cQ}\equiv 0$), such factorisation through the master\/-\/equation $\schouten{\cP_{0},\cP_{0}}=0$ can be immediate: here, we have that
%\footnote{Otherwise speaking, the flow $\frac{\Id}{\Id \varepsilon}(\cP)=\schouten{\cP,\cX(\cP)}$ is tautologically Poisson with respect to its native structure~$\cP$.} 
$%\[%\begin{align}
\schouten{\cP_{0},\cQ}=\schouten{\cP_{0},\schouten{\cP_{0},\cX}}=\tfrac{1}{2}\schouten{\schouten{\cP_{0},\cP_{0}},\cX}=\bigl(\tfrac{1}{2}\schouten{\cdot,\cX}\bigr)\,\bigl(\schouten{\cP_{0},\cP_{0}}\bigr)
$ %\]%\end{align}
for all $\boldsymbol{\partial}_{\cP_0}$-\/exact infinitesimal deformations $\cQ=\boldsymbol{\partial}_{\cP_{0}}(\cX)$ of the Poisson bi\/-\/vectors~$\cP_{0}$. Elaborated in~\cite{sqs15}, the Poisson cohomology estimate mechanism of the vanishing $\schouten{\cP_{0},\cQ}\doteq 0$ via $\schouten{\cP_{0},\cP_{0}}=0$ works --\,for the nontrivial cocycles $\cQ \notin \text{im}\,\boldsymbol{\partial}_{\cP_{0}}$ in the $\boldsymbol{\partial}_{\cP_{0}}$-\/cohomology\,-- due to much more refined principles. 
%We shall address %it 
%this mechanism in a subsequent paper.
That vanishing mechanism is applied to the factorisation problem at hand
in the paper~\cite{f16} (joint with R.\,Buring), where we prove the above conjecture.

\ack%\subsubsection*{Acknowledgements}
The second author thanks the Organizing committee of XXIV~International conference `Integrable systems \& quantum symmetries'
(13\,--\,19 June 2016; CVUT Prague, Czech Republic) for a warm atmosphere during the meeting. 
%%%
The research of~A.\,B. was supported by JBI~RUG project\:190.135.021;
A.\,K.~was supported by NWO grant VENI~639.031.623 (The Netherlands) and 
JBI~RUG project~106552 (Groningen%, The Netherlands
). 
A.\,B.~thanks R.\,Bu\-r\-ing for fruitful cooperation; %and communication; 
A.\,K.~thanks M.\,Kon\-tse\-vich for posing the problem and stimulating discussion.

\appendix
\section{The mechanism of vanishing for $\schouten{\cP,\cQ_{1:6}(\cP)}=0$: an example}\label{Perturbation theory}%
We wish to recognize the differential consequences of the Jacobi identity in the compatibility equation $\schouten{\cP,\cQ_{1:6}(\cP)}=0$, to understand why it holds. By a straightforward calculation we learn that $\schouten{\cP,\cQ_{1:6}(\cP)}=0$ for all Poisson bi-vectors on $\mathbb{R}^{3}$. But as soon as the differential consequences of the Jacobi identity are recognized, they can be translated into graphs. Independent of dimension, the language of graphs then answers the question which we started out with. This answer is found in~\cite{f16}.

Let us now illustrate a more analytic approach to the factorization problem for $\schouten{\cP,\cQ_{1:6}}=0$ via $\schouten{\cP,\cP}=0$ (see~\cite[App.~D]{f16} for details). The compatibility equation is a vanishing expression, which is impossible to factorize through the Jacobi identity, which itself is also zero. To make both visible, we perturb a given Poisson bi-vector $\cP$ using $\tilde{\cP}=\cP+\epsilon \cdot \Delta$ for a bi-vector $\Delta$, in such a way that $\tilde{\cP}$ is no longer Poisson, thereby $\schouten{\tilde{\cP},\tilde{\cP}} \neq 0$. The goal is to perturb the bi-vector $\cP$ such that the left-hand side  $\schouten{\tilde{\cP},\tilde{\cQ}_{1:6}}$ becomes non-zero as well. Now the Jacobi identity's non-zero differential consequences becomes recognizable in the non-zero expression $\schouten{\tilde{\cP},\tilde{\cQ}_{1:6}}$.

\begin{example}
Consider the Poisson bi-vector obtained on $\mathbb{R}^{3}$ from the determinant construction using two functions $g(z)$ and $f(x)$ as argument and pre-multiplication factor, respectively. Let the perturbation $\Delta$ be given component-wise by $\Delta^{12}=f_{1}(y,z)$, $\Delta^{13}=f_{2}(y,z)$ and $\Delta^{23}=0$.
The perturbed bi-vector then equals
\[
\tilde{\cP}=\left[\begin{matrix}
0 & f\cdot \Id g/\Id z %\frac{\partial g}{\partial z} 
& 0\\
-f\cdot \Id g/\Id z %\frac{\partial g}{\partial z} 
& 0 & 0\\
0 & 0 & 0
\end{matrix}\right]
+
\epsilon \cdot \left[\begin{matrix}
0 &  f_{1}& f_{2}\\
-f_{1} & 0 & 0\\
-f_{2} & 0 & 0
\end{matrix}\right].
\]
The left\/-\/hand sides of the Jacobi identity and of the compatibility condition are evaluated to
\[%\begin{eqnarray*}
\schouten{\tilde{\cP},\tilde{\cP}}^{123}=%&=&
\epsilon f_{2}\cdot%\left(
\frac{\Id f}{\Id x}%\cdot\right)\left(
\frac{\Id g}{\Id z}%\right)
+\bar{o}(\epsilon),
\qquad%\\
\schouten{\tilde{\cP},\tilde{\cQ}}^{123}=%&=&
-\epsilon\cdot%\left(
\frac{\partial^{3}f_{2}}{\partial y^{3}}%\right)
\left(\frac{\Id f}{\Id x}\right)^{4}\left(\frac{\Id g}{\Id z}\right)^{4}+\bar{o}(\epsilon).
\]%\end{eqnarray*}
There is only one way to recognize a differential consequence of the Jacobiator inside $\schouten{\tilde{\cP},\tilde{\cQ}_{1:6}}^{123}$. Namely, the Jacobi identity contains a product of $f_{2}$ and derivatives of $f$ and $g$. The same is true for its non-zero differential consequences. Let us extract this product from $\schouten{\tilde{P},\tilde{Q}_{1:6}}^{123}$. The only differential consequences of $f_{2}$, $\Id f/\Id x$, and $\Id g/\Id y$ in $\schouten{\cP,\cQ_{1:6}}^{123}$ are $\partial^{3} f_{2} / \partial y^{3}$, $\Id f/\Id x$ and $\Id g/\Id z$, respectively. This hints that we have the differential consequence $\schouten{\cP,\cP}^{123}_{yyy}$. To understand what its coefficient is, we note that the remaining co\/-\/factors in $\schouten{\tilde{\cP},\tilde{\cQ}_{1:6}}^{123}$ form %equal 
$\smash{(\cP^{12}_{x})^{3}}$. We conclude that the left\/-\/hand side of the compatibility equation factorizes through the Jacobi identity as follows
\[\schouten{\cP,\cQ_{1:6}}^{123} = \cP^{12}_{x}\cP^{12}_{x}\cP^{12}_{x}
\lshad\cP,\cP\rshad^{123}_{yyy}+ \cdots.
\]
Looking at this expression, we construct a list of graphs that can encode it. Such a list %once 
fully formed, it is subtracted from $\schouten{\cP,\cQ_{1:6}}$ and resolved with respect to the coefficients of every proposed graph. We keep subtracting the already found graphs from any non-zero perturbations of $\schouten{\cP,\cQ_{1:6}}$ in the future, once the coefficients are known. The example under study gave us the tripod graph, which is the first entry 
%in equation (6) 
in~\cite[Eq.\:(6)]{f16}. Proceeding in the same way, we also recognized the 'elephant' graph, which is the sixth entry in that solution (cf.~\cite[Remarks 10--11]{f16}).
\end{example}

%The author thanks the Organising committee of the 7th %seventh 
%International workshop
%`Group Analysis of Differential Equations and Integrable Systems' 
%(15--19 June 2014, Larnaca, Cyprus) for %hospitality.
%a warm atmosphere during the meeting.
%
%The author also thanks the students who followed the master course `Differential calculus in (non)\/commutative geometry' (JBI RUG, 2013) 
%for stimulating discussions and constructive criticism.

%\medskip
%\bibliographystyle{iopart-num}
%\bibliography{cyprus}

\section*{References}
\begin{thebibliography}{10}
\expandafter\ifx\csname url\endcsname\relax
  \def\url#1{{\tt #1}}\fi
\expandafter\ifx\csname urlprefix\endcsname\relax\def\urlprefix{URL }\fi
\providecommand{\eprint}[2][]{\url{#2}}
% Bibliography created with iopart-num v2.0
% /biblio/bibtex/contrib/iopart-num

\bibitem{Ascona96}
{Kontsevich M} (1997)
Formality conjecture
\textit{Deformation theory and symplectic geometry (Ascona %, June 17--21 
1996)}
%\textit{Math.\ Phys.\ Stud.}~
\textbf{20} (Dordrecht: Kluwer Acad.\ Publ.) 139--156

\bibitem{KontsevichFormality}
{Kontsevich M} (2003)
{Deformation quantization of {P}oisson manifolds}
\textit{Lett.\ Math.\ Phys.} \textbf{66}:3 157--216
(\textit{Preprint} \texttt{q-alg/9709040})

\bibitem{MKZurichICM}
{Kontsevich M} (1995)
Homological algebra of mirror symmetry
\textit{Proc.\ Intern. Congr. Math.} \textbf{1} %,\,\vol{2} 
(Basel: Birkh\"auser) 120--139

\bibitem{MKParisECM}%{Kontsevich1994}
{Kontsevich M} (1994)
Feynman diagrams and low\/-\/dimensional topology
\textit{First Europ.\ Congr.\ of Math.} \textbf{2} 
(Paris, 1992)
%{Progr.\ Math. 
\textbf{120} (Basel: Birkh\"auser) 97--121
%(Reviewer: Anatoly Libgober) 57R57 (14H15 32G15 57M25)

\bibitem{sqs15}
{Buring R and Kiselev A V} (2017) On the Kon\-tse\-vich $\star$-\/product asso\-ci\-a\-ti\-vi\-ty me\-cha\-n\-ism\ 
\textit{PEPAN Letters} \textbf{14}:2 accepted %4\,pp.\ 
(\textit{Preprint} \texttt{arXiv:1602.09036} [q-alg])

\bibitem{f16}
{Bouisaghouane A, Buring R and Kiselev A V} (2016)
The Kontsevich tetrahedral flows revisited
\textit{Preprint} \texttt{arXiv:1608.01710} (v2) [q-alg] 20%~p.

\bibitem{Donin}
{Donin J} (1998)
On the quantization of quadratic Poisson brackets on a polynomial algebra of four variables
\textit{Lie Groups and Lie Algebras. Their representations, generalisations and applications}
%(B.\,P.\,Kom\-ra\-kov, I.\,S.\,Kra\-sil'\-shchik, G.\,L.\,Lit\-vi\-nov and A.\,B.\,Sossin\-sky, eds), 
%Math. Appl. 
\textbf{433}
(Dordrecht: Kluwer Acad.\ Publ.) %Springer Sci.\ Business Media, 
17--25

\bibitem{Perelomov}
{Grabowski J, Marmo G and Perelomov A M} (1993)
Poisson structures: towards a classification
\textit{Mod.\ %ern 
Phys.\ Lett.} \textbf{A8}:18  1719--1733

%\bibitem{Nutku}
%{G\"{u}mral H, Nutku Y} (1993) Poisson structure of dynamical systems with three degrees of freedom, 
%{J.~Math.\ Phys.}~vol 34:12, 5691--5723.

\bibitem{Van}
{Vanhaecke P} (1996) 
\textit{Integrable systems in the realm of algebraic geometry} 
%Lect.\ Notes Math.\ 
\textbf{1638} %(2nd ed.), 
(Berlin: Springer\/--\/Verlag) %Heidelberg. %xxviii+
%513

\bibitem{Gerstenhaber}
{Gerstenhaber M} (1964) 
On the deformation of rings and algebras
\textit{Ann.\ Math.} \textbf{79} 59--103

\bibitem{Olver}
{Olver P J} (1993) 
\textit{Applications of Lie groups to differential equations}
%{Grad.\ Texts in Math.} 
\textbf{107} (2nd ed.) 
(NY: Springer\/--\/Verlag) %xxviii+
%513

\bibitem{Kent}%https://kar.kent.ac.uk/8141/1/Wang_List_1%2B1.pdf
{Wang J P} (2002)
A list of $1+1$ dimensional integrable equations and their properties
\textit{J.~Nonlin.\ Math.\ Phys.} \textbf{9} %suppl.~1
%Recent advances in integrable systems (Kowloon, 2000), 
213--233

\bibitem{gvbv}
{Kiselev A V} (2013)
The geometry of variations in {B}atalin\/--\/{V}ilkovisky formalism
\textit{J.~Phys.: Conf.\ Ser.} {\bf 474} Paper~012024 1--51\ %
(\textit{Preprint} \texttt{1312.1262} \texttt{[math-ph]})

\bibitem{dq15}
{Kiselev A V} (2015) 
Deformation approach to quantisation of field models
\textit{Preprint} $\smash{\text{IH\'ES}}$/M/15/13 (Bures\/-\/sur\/-\/Yvette, %2015
France) 37

\bibitem{prg15}
{Kiselev A V} (2016)
The right\/-\/hand side of the Jacobi identity: to be naught or not to be\,?
\textit{J.~Phys.\textup{:}\ Conf.\ Ser.} \textbf{670} 
   %Proc.\ XXIII Int.\ conf.\
   %`Integrable Systems and Quantum Symmetries' (23--27 June 2015, CVUT Prague,
   %Czech Republic), 
%Paper~
012030  1--17 (\textit{Preprint} \texttt{arXiv:1410.0173} [math-ph])

\bibitem{Vinogradov1997}
{Vinogradov A and Vinogradov M} (1998)
On multiple generalizations of Lie algebras and Poisson manifolds
\textit{Secondary calculus and cohomological physics} 
%(Moscow 1997,M.~Henneaux, J.~Krasil'shchik and A.~Vinogradov, eds),
%Contemp.\ Math.\ 
\textbf{219} (Providence RI: AMS) 273--287

\bibitem{Nambu}
{Nambu Y} (1973) 
Generalized Hamiltonian dynamics 
\textit{Phys.\ Rev.~D} \textbf{7} 
%American Physical Society, 
2405--2412

\bibitem{Yoshioka}
{Omori H, %Hideki(J-SUTST)
Maeda Y  %Yoshiaki(J-KEIOES); 
and Yoshioka A} %Akira(J-SUTST)
(1993)
A construction of a deformation quantization of a Poisson algebra
\textit{Geometry and its applications} 
%(Yokohama 1991,
%Edited by 
%T~%Tadashi 
%Nagano, H~%Hideki 
%Omori, Y~%Yoshiaki 
%Maeda and M~%asahiko 
%Kanai eds), 
(River Edge NJ: World Sci.\ Publ.) 201--218

\bibitem{cycle14}
{Kiselev A V} (2015)
The calculus of multivectors on noncommutative jet spaces 
\textit{Preprint} $\smash{\text{IH\'ES}}$/M/14/39 (Bures\/-\/sur\/-\/Yvette, %2014, 
France) \texttt{arXiv:1210.0726} (v3) [math.DG] 41

\bibitem{Olver1997}
{Fokas A S, Olver P J and Rosenau P} (1997)
A plethora of integrable bi\/-\/Hamiltonian equations
\textit{Algebraic aspects of integrable systems}
% In Memory of Irene Dorfman (
%\textit{Progr.\ Nonlin.\ Diff.\ Equations and Their Appl.\ } 
\textbf{26}
%A S Fokas and I M Gelfand eds, 
(Boston MA: Birkh\"{a}user) %Boston, 
93--101

\bibitem{Vodova}
{Vodov\'a J} (2013)
Low\/-\/order Hamiltonian operators having momentum
\textit{J.~Math.\ Anal.\ Appl.} \textbf{401}:2 724--732\ %
(\textit{Preprint} \texttt{arXiv:1111.6434} [math-ph])

\bibitem{Merkulov}
{Merkulov S A} (2010)
Exotic automorphisms of the Schouten algebra of polyvector fields
\textit{Preprint} \texttt{arXiv:0809.2385} (v6) [q-alg] 37

\end{thebibliography}
%\end{document}

\providecommand{\newblock}{}

\end{document}